\documentclass[12pt]{article}
\usepackage{amsmath,amssymb,latexsym,graphicx}
\newtheorem{theorem}{Theorem}[section]
\newtheorem{lemma}[theorem]{Lemma}

\newtheorem{remark}[theorem]{Remark}
\newtheorem{corollary}[theorem]{Corollary}
\newtheorem{definition}[theorem]{Definition}
\newtheorem{conjecture}[theorem]{Conjecture}

\def\pf{{\bf Proof }}
\numberwithin{equation}{section}

\begin{document}
\title{Locally polynomially integrable surfaces and finite stationary phase expansions}
\author{Mark ~Agranovsky}
\maketitle
\begin{abstract}
Let $M$ be a strictly convex smooth connected hypersurface in $\mathbb R^n$ and $\widehat{M}$ its convex hull. We say that $M$ is polynomially integrable   if the $(n-1)-$ dimensional volumes of the sections of $\widehat M$ by hyperplanes, sufficiently close to the tangent hyperplanes to $M,$ depend polynomially on the distance of the hyperplanes to the origin. It is conjectured that only quadrics in odd dimensional spaces possess such a property. The main result of this article partially confirms the conjecture. The study of integrable domains and surfaces is motivated by  a conjecture of V.I. Arnold about algebraically integrable domains. The result and the proof are related to study oscillating integrals for which the asymptotic stationary phase expansions consist of finite number of terms.

\end{abstract}

\section{Introduction}\label{S:Intro}
In \cite{Ag}, the following definition was introduced for bounded domains in $\mathbb R^n.$
The body $K \subset \mathbb R^n$  is called {\it polynomially integrable} if the Radon transform of its characteristic function $\chi_K:$
$$A_K(\xi,p)=vol_{n-1}\big( \{ x \in K: \langle x, \xi \rangle =p \} \big)= \int\limits_{ \langle x, \xi \rangle=p  } \chi_K(x) dV_{n-1}(x)$$
is a polynomial in $p.$  The same term will be addressed the boundary $\partial K$ which in our case will be assumed  a smooth hypersurface.

The function $A(\xi,p),$ which we will call {\it the (sectional) volume function,} evaluates the $(n-1)-$ dimensional volume of the cross-section of $K$
by the hyperplane $\{x \in \mathbb R^n: \langle x, \xi \rangle=p \}.$

It was proved that there is no polynomially integrable bodies with $C^{\infty}$ boundary  in $R^{2k},$ while in odd dimensions only solid ellipsoids are polynomially \cite{KMY}, \cite{Ag}, and even rationally or real-analytically \cite{Ag1}, integrable domains with $C^{\infty}$ boundary.

In this article, we extend the notion of polynomial integrability to the case of open hypersurfaces. In this case, the  volume
function $A_M(\xi,p)$ is defined for the cross-section of $M$ by hyperplanes close to tangent those. To guarantee the finiteness of the volumes, we assume that $M$ is strictly convex.

The main question that we are concerned with is: {\it what are polynomially integrable hypersurfaces?}
The quadrics in odd dimensional spaces deliver examples of such hypersurfaces and it is expected that there are no other examples.

The main result of this article is towards this conjecture. Namely,
we prove that, under certain additional conditions, the polynomiality of the volume function $A(\xi,t)$ implies that the hypersurface $M$ is a quadric, namely, an elliptic paraboloid, in $\mathbb R^{2m+1}.$

The study of polynomially integrable bodies and surfaces has been motivated by works \cite{Ar},\cite{Vas},  devoted to algebraically integrable bodies. They are bodies $K$  for which the two-valued volume function $V^{\pm}_K(\xi,t),$ where $ V^+(\xi,t)=\int_{-\infty}^{t} A(\xi,u)du, V^{-}(\xi,t)=\int_t^{+\infty}A(\xi,u)du,$  is algebraic.
Celebrated Newton's Lemma about ovals from \cite{N} states that there is no such domains (with infinitely smooth boundary) in the plane.
V.I. Arnold (\cite{Ar},problems 1987-14, 1988-13, 1990-27) has suggested to generalize Newton's lemma to higher dimensions and conjectured that there is no algebraically integrable domains with smooth boundary in even-dimensional spaces, while in odd-dimensional spaces the only such domains are ellipsoids. The first part of the conjecture was confirmed by Vassiliev  \cite{Vas1,Vas,VasBook}, where he has proved that there is no algebraically integrable domain with smooth boundary  in $\mathbb R^n$ when $n$ is even. The case of odd dimensions is not completely solved yet.
The main result of this article can be viewed as one in the direction of the Arnold's conjecture.

\section{Definitions and main result}\label{S:Main}

Let $M$ be a differentiable strictly convex open connected hypersurface in $\mathbb R^n.$

Denote $\widehat M$ the convex hull of $M.$ Let $a \in M,$ $T_a(M)-$ the tangent hyperplane at $ M,$ and  $\nu_a-$ the unit normal vector to $M$ directed towards $\widehat M.$

When $t$ is positive and small then the boundary of the intersection of $\widehat M$ with the translate of the tangent plane is contained in $M:$
$$\widehat M \cap (T_a(M)+t \nu_a) \subset \widehat M$$
and the $(n-1)-$ dimensional volume of the intersection
\begin{equation}\label{E:Aa}
A_M^a (t)=vol_{n-1} \big(\widehat M \cap (T_a(M)+t\nu_a)\big)=\int\limits_{ \{x \in \widehat M: \langle x-a,\nu_a \rangle =t \} } dV_{n-1}.
\end{equation}
is finite.

\begin{remark}
Notice that we have re-parameterized the volume function $A(\xi,p)$ defined in Introduction by introducing the parameter $t=p-\langle \xi,a \rangle$ instead of $p.$ The parameter $t$ expresses the distance of the hyperplane to the tangent one, while $p$ is the distance of the hyperplane to the origin. Obviously, the properties of polynomiality with respect to either parameter are equivalent.
Everywhere in the sequel we use the parameter $t.$
\end{remark}

\begin{definition}
The hypersurface $M$ is called (locally) polynomially integrable if for any $a \in M$ there exists $\varepsilon_a >0$ such that the volume function $ A_M^a (t)$ is a polynomial for $ t \in [0,\varepsilon_a).$
\end{definition}
\begin{figure}[h]
\scalebox{0.5}{\includegraphics{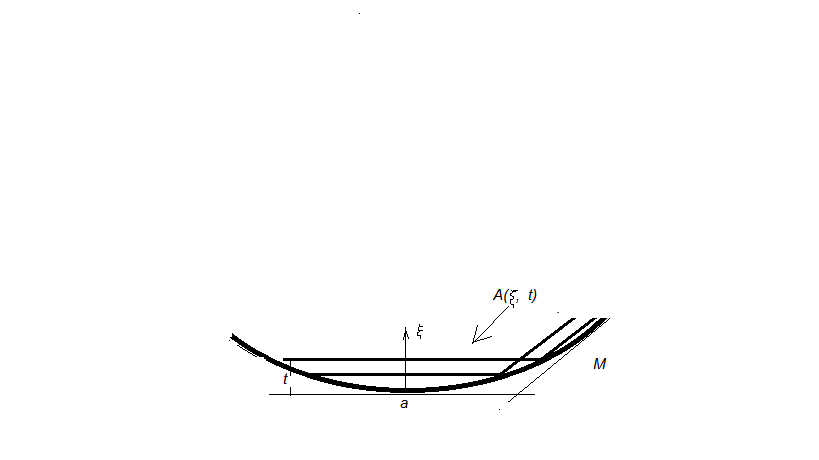}}
\label{Fig.1}
\end{figure}
It is not difficult to check that there is  no such hypersurfaces in $\mathbb R^{2m}.$
Indeed, $M$ is strictly convex and hence there is a non-degenerate (elliptic)  point $a \in M$ at which $M$ has strictly positive Gaussian curvature $\kappa >0.$ Then
\begin{equation} \label{E:kappa}
A_M^a(t)=\kappa t^{\frac{n-1}{2}}(1+o(1)), t \to 0+,
\end{equation}
\cite{GG},\cite {Ag}.
Since for even $n$ the exponent $\frac{n-1}{2}$ is non-integer, the function $A_M^a(t)$ cannot be a polynomial near $t=0.$

However, there are polynomially integrable strictly convex hypersurfaces in odd dimensions. They are convex connected quadrics in $\mathbb R^{2m+1}$ which can be transformed by translations and rotations to one of the following surfaces:
\begin{enumerate}
\item ellipsoid  $a_1^2 x_1^2+ \cdots a_n^2 x_n^2=1,$
\item two sheet hyperboloid $a_1^2 x_1^2+\cdots +a_{n-1}^2 x_{n-1}^2- a_n^2 x_n^2= -1, x_n > 0,$
\item elliptic paraboloid $a_1^2 x_1^2+\cdots a_{n-1}^2 x_{n-1}^2 - a_n^2 x_n=0.$
\end{enumerate}
The polynomial integrability of each above quadratic hypersurface can be checked by a straightforward computation.

\begin{conjecture} \label{C:conj}  The only locally integrable strictly convex $C^{\infty}$  hypersurfaces are, up to an affine transformation, the above enlisted  quadrics in odd-dimensional spaces.
\end{conjecture}

In this article, we  characterize the surfaces of the third type, i.e. elliptic paraboloids, in terms of polynomial integrability.
To formulate the result, we need to remind the notion of osculating paraboloid. Everywhere, the hypersurface $M$ is assumed at least $C^{\infty}.$

Let $ a \in M.$ Represent, near $a,$ the hypersurface $M$ as the graph, say,
$$x_{n}=f(x^{\prime})$$
where
$$ x^{\prime}=(x_1, \dots,x_{n-1}).$$
and $f$ is a $C^{\infty}$ function defined in a neighborhood of $a^{\prime}.$

The osculating paraboloid $P_a$ is defined as the graph of the second order Taylor polynomial  of the function $f(x^{\prime})$ at the point $a^{\prime}:$
$$P_a: x_{n}=f(a)+df_{a^{\prime}}(x^{\prime}-a^{\prime})+\frac{1}{2}d^2 f_{a^{\prime}}(x^{\prime}-a^{\prime}).$$

The point $a \in M$ is {\it non-degenerate} if the second differential $d^2f_{a^{\prime}}$ is a non-degenerate quadratic form, i.e.,
is Hessian at $a^{\prime}$ is different from 0. In this case, the Gaussian curvature at $a$ is not $0$.
If $M$ is strictly convex in a neighborhood of a non-degenerate point $a$  then all the principal curvatures at $a$ are non-zero and of the same sign. Let $k \in \mathbb N, k > 2.$ We say that the osculating paraboloid $P_a$ has {\it contact with $M$ of order higher than $k$} if
$$df^3_{a^{\prime}}=... = df^k_{a^{\prime}}=0.$$

In this article, we  prove that following result towards Conjecture \ref{C:conj}:
\begin{theorem}\label{T:main}
Let $M$ be an open real analytic strictly convex polynomially integrable hypersurface in $\mathbb R^n, n \geq 3$ is odd.
Suppose that there exists a non-degenerate point $a^0 \in M$ at which the osculating paraboloid contacts $M$ to the order higher than four. Then $M$ is an elliptic paraboloid, i.e. $M$ can be transformed , by a suitable affine transformation,to the surface
$$x_n=\sum_{j=1}^{n-1} \alpha_j^2 x_j^2.$$
\end{theorem}

\subsection{Plan of the proof}
First, we prove that the condition of polynomiality of the volume function $A_M^a$ is equivalent to the finiteness of the stationary phase expansion of certain oscillating Fourier integrals on $M.$  In order to get rid of the contribution of the boundary, we use a family of cut-off functions. Passing to the normal Morse coordinates on $M$ near the distinguished point $a^0,$ we make the phase function quadratic. Then we apply to the stationary phase asymptotic expansions the iterates of the wave operator $\Box_{\omega}$  with respect to the variable $\omega,$ which expresses the direction of the normal vector $\nu_a,$ and  then evaluate the result at the point $a^0.$  We observe that, due to the osculating paraboloid condition  at $a^0,$  when we repeatedly apply the operator $\Box_{\omega}$ to the oscilat8ng integral, then the power of the reciprocal large parameter in the leading term of the stationary phase expansion grows faster than that of the last term. Therefore, after applying $\Box_{\omega}$ finite number of times, the finite part of the expansion disappears.
That translates as vanishing identically on $M$  the symbol $x_n-\sum_{j=1}^{n-1}x_j^2$ of the wave operator $\Box_{\omega}$ which implies that $M$ satisfies the equation of an elliptic paraboloid.

\section{Stationary phase expansion}

In this section we will relate the property of polynomial integrability with a property of finiteness of stationary phase expansion of a Fourier integral. Everywhere in the sequel, the dimension $n$ of the ambient space is assumed odd.

Denote
$$\gamma:M \to S^{n-1}$$
the Gaussian mapping that maps a point in $M$ to the unit, inward with respect to $\widehat M,$ normal vector $\nu_a$ at that point:
$$\gamma(a)=\nu_a.$$

Theorem  \ref{T:main} is local, because $M$ is assumed real analytic and connected, hence we can reduce $M$ to a neighborhood of the distinguished point $a^0$ mentioned in the formulation of Theorem \ref{T:main}. The hypersurface $M$ is strictly convex and $a^0$ is a non-degenerated point, hence the neighborhood can be chosen so that the Gaussian mapping  $\gamma$ is a real analytic diffeomorphism
$$\gamma: M \to \gamma(M) \subset S^{n-1}$$ of $M$ onto an open subset of the unit sphere.

It will be convenient to parametrize the volume function $A_M^a(t),$ defined in (\ref{E:Aa}), by the normal vector $\gamma(a)$ and introduce
\begin{equation} \label{E:AMxi}
A_M(\xi,t)=A_M^{\gamma^{-1}(\xi)}(t).
\end{equation}

Introduce the function
\begin{equation}\label{E:h_M}
h_M(\xi)=\langle \gamma^{-1}(\xi),\xi \rangle.
\end{equation}.
The absolute value of this function evaluates the distance to the origin of the tangent plane with the normal vector $\xi:$ $|h_M(\xi)|=dist(T_a(M),0), \nu_a=\xi.$

Since $M$ is strictly convex, for all $x \in M, x \neq a,$ holds
$$\langle x-a, \xi \rangle >0$$ and hence
$$\langle x,\xi \rangle > h_M(\xi)$$
with the equality for $x=a.$ Here $\xi=\gamma(a)$ is the unit inward normal vector at $a.$

\subsection{ Cut-off functions $\rho_c$}
Fix $c>0.$ Let $\rho_c(t)$ be a function with the following properties:

\begin{enumerate}
\item $\rho_c \in C^{\infty}(\mathbb R).$
\item $ \rho_c(t)=1$ in a neighborhood of $t=0.$
\item $ supp \rho_c  \subset (-c,c). $
\end{enumerate}

Now we define the function $\rho_{c,\xi}(x), x \in \mathbb R^n$ by:
\begin{equation} \label{E:rho}
\rho_{\xi,c}(x):=\rho_c \big(\langle x-a, \xi  \rangle )=\rho_c(\langle x,\xi \rangle - h_M(\xi)\big),
\end{equation}
where $a \in M $ is the point with the normal vector $\xi,$ i.e., $a=\gamma^{-1}(\xi),$ and  $h_M(\xi)$ is the "support" function (\ref{E:h_M}),
$h_M(\xi)=\langle a, \xi\rangle .$  The function $\rho_{\xi,c}(x)$ is supported in the strip $ | \langle \xi,x\rangle - h_M(\xi)| < c.$

Define the cap $M_{\xi,c} \subset M$ by
\begin{equation}\label{E:M^a_c}
M_{\xi,c}=\{ x \in M:  \langle x, \xi \rangle-h_M(\xi)  \leq c \}.
\end{equation}

We assume that $c > 0$ is so small that the set $M_{\xi,c}$
is relatively compact in $M.$ Then the restriction of $\rho_{c,\xi}$ to $M$ is supported strictly inside $M.$

Let us summarize the properties of the family of cut-off functions $\rho_{c,\xi}:$
\begin{enumerate}
\item $\rho_{\xi,c} \in C^{\infty}(\mathbb R^n),$
\item  $\rho_{\xi,c}(x)=1$ for $x$ in a neighborhood of $a=\gamma^{-1}(\xi),$
\item $M \cap supp \rho_{\xi,c} \subset M_{\xi,c},$
\item $M_{\xi,c}$ is relatively compact in $M.$
\end{enumerate}

\subsection{Oscillating integral $I_{\xi,c}$}

Consider the integral
\begin{equation}\label{E:I}
I_{\xi,c}(\lambda):=\int\limits_M \frac{\partial}{\partial \nu_x} e^{i \lambda  \langle \xi , x \rangle }
\rho_{\xi,c}(x) dS(x),
\end{equation}
where $dS(x)$ is the surface measure on $M$ and $\nu_x$ is the unit normal vector to $M$ at $x$ directed to the convex hull $\widehat M.$

The integral $I_{\xi,c}$ can be written as
$$
I_{\xi,c}(\lambda) =\int\limits_M e^{i \lambda  \langle \xi , x \rangle } i\lambda \langle \xi,\nu_x \rangle \rho_{\xi,c}(x) dS(x)
$$
and therefore can be viewed as an oscillating integral, with respect to the large parameter $\lambda$ and with the phase function
$$\varphi_{\xi}( x)=\langle \xi, x \rangle. $$

Notice that due to the location of the support of the cut-off function $\rho_{c,\xi},$ the integration in (\ref{E:I}) is performed, in fact, over the relatively compact submanifold $M_{\xi,c} \subset M,$  defined in (\ref{E:M^a_c}).

According to the stationary phase method, the asymptotic of $I_{\xi,c}(\lambda),$ as $\lambda \to \infty,$ is determined by the  values of the weight function
$$
M_{\xi,c} \ni x \to \langle \xi,\nu_x \rangle \rho_{\xi,c}(x),
$$
near critical points of the phase function $\varphi_{\xi}$ on $M_{\xi,c},$ i.e., point $x \in M_{\xi,c}$ such that $\nu_x=\xi.$
In our case,  the only such point is
$$a=\gamma^{-1}(\xi)$$
The value of the phase $\varphi_{\xi}$ function at its critical point $a=\gamma^{-1}(\xi)$ is
$$\varphi_{\xi}(a)= \langle a, \nu_a \rangle=h_M(\xi),$$
where $h_M$ is the "support" function (\ref{E:h_M}).

Then the stationary phase method (see e.g. \cite{Wong}, Ch.IX, Thm.1) yields the following expansion of the oscillating integral $I_{\xi,c}(\lambda)$ with the large parameter $\lambda:$
\begin{equation}\label{E:expansion}
I_{\xi,c}(\lambda)=const \ e^{i \lambda h_M(\xi)}\frac{1}{ \lambda^{\frac{n-1}{2}}} \sum_{j=0}^{\infty} \frac{b_k(\xi)}{\lambda^{j}},
\end{equation}

The series in (\ref{E:expansion}) an is asymptotic series, i.e.,  for any $N$
$$I_{\xi,c} (\lambda)=const \  e^{i h_M(\xi)} \lambda^{-\frac{n-1}{2}} \Big( S_N(\xi,\lambda)+ o\big( \frac{1}{\lambda^N} \Big) , \lambda \to \infty,$$
where $S_N(\xi,\lambda)$ is the $N-th$ partial sum of the formal series in (\ref{E:expansion}).

\subsection{Finiteness of the stationary phase expansion of $I_{\xi,c}(\lambda)$}

\begin{lemma} \label{L:oscill} If $M$ is  polynomially integrable then the asymptotic series (\ref{E:expansion})
contains only finite number of nonzero terms, i.e., there exists $\overline N$ such that $S_N=S_{\overline{N}}$ for all $ N \geq \overline{N}$ and hence
the above decomposition takes the form
\begin{equation} \label{E:decomp}
I_{\xi,c}(\lambda)=const e^{ih_M(\xi)} Q_{\xi} \Big(\frac{1}{\lambda} \Big)+R_{\xi}(\lambda),
\end{equation}
where $Q_{\xi}$ is a polynomial of $\deg Q_{\xi} \leq N_0=\overline {N} +\frac{n-1}{2}$ with coefficients depending on $\xi,$ and for any natural $N$ holds $R_{\xi}(\lambda)=o \Big(\frac{1}{\lambda^N}\Big), \lambda \to \infty,$
uniformly, with all $\xi-$ derivatives,  with respect to the parameter $\xi.$
\end{lemma}
\pf
Denote $D_{\xi,c}$  the portion of the convex hull $\widehat M$ cut off by the hyperplane parallel to the tangent hyperplane to $M,$ with the normal vector $\xi$ and the distance $c>0$ from the tangent hyperplane:
$$D_{\xi,c}= \{x \in \widehat M:  \langle \xi, x-a \rangle =\langle \xi,x \rangle -h_M(\xi) \leq  c \},$$
where $a=\gamma^{-1}(\xi)$ is the point in $M$ with the unit normal vector $\nu_a=\xi.$

The  parameter $c$ is chosen so small that
$D_{\xi,c}$ is compactly contained in $\widehat M$. Then the boundary of the domain $D_{\xi,c}$ is the union of the two parts:
$$
\partial D_{\xi,c}=M_{\xi,c} \cup \Omega_{\xi,c},
$$
where $M_{\xi,c}$ is the cap defined in (\ref{E:M^a_c} ) and $ \Omega_{\xi,c}$ is the "top cover"
$$\Omega_{\xi,c}=\{ x \in \widehat M: \langle \xi, x-a \rangle = \langle \xi,x \rangle -h_M(\xi)= c \}, a=\gamma^{-1}(\xi).$$

Now apply the Stokes formula to the volume integral over the domain $D_{\xi,c}:$
$$
\begin{aligned}
&\int\limits_{D_{\xi,c}}  \Delta e^{i\lambda \langle \xi, x \rangle}\rho_{\xi,c}(x)dV(x)
-\int\limits_{D_{\xi,c}}  e^{i\lambda \langle \xi, x \rangle} \Delta \rho_{\xi,c}(x) dV(x) \\
&=\int\limits_{\partial D_{\xi,c}} \frac{\partial}{\partial \nu_x} e^{i\lambda \langle \xi, x \rangle}\rho_{\xi,c}(x)dS(y)
-\int\limits_{\partial D_{\xi,c}}  e^{i\lambda \langle \xi, x \rangle} \frac{\partial}{\partial \nu_x} \rho_{\xi,c}(x)dS(x).
\end{aligned}
$$

Compute Laplace operator of the exponential factor in the first integral in the left hand side and express the first surface integral in the right hand side through the other three integrals participating in the identity. One obtains:
\begin{equation}\label{E:F123}
\int\limits_{\partial D_{\xi,c}} \frac{\partial}{\partial \nu_x}  e^{i\lambda \langle \xi, x \rangle}\rho_{\xi,c}(x)dS(x)
=F_{\xi,1}(\lambda)+F_{\xi,2}(\lambda)+F_{\xi,3}(\lambda),
\end{equation}
where
\begin{equation}
\begin{aligned}
&F_{\xi,1}(\lambda)=-\lambda^2 \int\limits_{D_{\xi,c}}  e^{i\lambda \langle \xi, x \rangle} \rho_{\xi,c} (x) dV(x),\\
&F_{\xi,2}(\lambda)=- \int\limits_{D_{\xi,c}}  e^{i\lambda \langle \xi, x \rangle} \Delta \rho_{\xi,c}(x) dV(x),\\
&F_{\xi,3}(\lambda)= \int\limits_{\partial D_{\xi,c}}  e^{i\lambda \langle \xi, x \rangle}  \frac{\partial}{\partial \nu_x} \rho_{\xi,c} (x) dS(x). \\
\end{aligned}
\end{equation}

The surface $\partial D_{\xi,c}$ of integration in the left hand side of (\ref{E:F123}) consists of two parts,
$M_{\xi,c}$  and $\Omega_{\xi,c}.$
However, $\Omega_{\xi,c}$ does not contribute to  integral (\ref{E:F123}), because the factor $\rho_{\xi,c}$ vanishes on the "cover" $\Omega_{\xi,c}.$
Therefore, the surface $\partial D_{\xi,c}$ of integration   in the left hand side of (\ref{E:F123}) can be replaced by $M_{\xi,c},$ and hence  (\ref{E:F123}) can be rewritten, due to definition (\ref{E:I}) of $I_{\xi,c}$, as
\begin{equation}\label{E:IF123}
I_{\xi,c}(\lambda)=F_{\xi,1}(\lambda)+F_{\xi,2}(\lambda)+F_{\xi,3}(\lambda).
\end{equation}

Now we start exploring the asymptotic expansions of each integral $F_{\xi,1}, F_{\xi,2}, F_{\xi,3}$ separately.

\subsubsection{ Asymptotic expansion of $F_{\xi,1}(\lambda)$}  \label{S:Fxi1}

By Fubini theorem, the volume integral $F_{\xi,1}(\lambda)$ can be reduced to a one-dimensional integrals:
$$
\begin{aligned}
&F_{\xi,1}(\lambda)=-\lambda^2  \int\limits_0^c \big( \int\limits_{\langle \xi, y \rangle -h_M(\xi)= t} e^{i\lambda \langle \xi, y \rangle } \rho_{\xi,c}(y) dy \big) dt \\
&=-\lambda^2 e^{i\lambda h_M(\xi)} \int_0^c e^{i\lambda t} A_M (\xi,t) \rho_c (t)dt,
\end{aligned}
$$
where $A_M (\xi,t)$ is the section volume function defined in (\ref{E:Aa}), (\ref{E:AMxi}).
In the computation, we have used also  the definition $\rho_{\xi,c}=\rho_c(t), t=\langle \xi, y \rangle -h_M(\xi).$

Thus, $F_{\xi,1}(\lambda)$ boils down to a one-dimensional Fourier oscillating integral over the segment $[0,c].$
Integration $N$ times by parts  yields:
\begin{equation}\label{E:sum}
\begin{aligned}
&-\lambda^2 \int_0^c e^{i\lambda t} \rho_c (t) A_M (\xi,t)dt \\
&=\sum_{k=0}^N (i\lambda)^{-k+1} (-1)^k e^{i \lambda t}
\frac{d^k}{dt^k} \big(\rho_c(t)A_M (\xi,t)\big) \big|^{t=c}_{t=0} +o(\lambda^{-N}).
\end{aligned}
\end{equation}
It follows from (\ref{E:kappa}) that the  index of summation in (\ref{E:sum})
satisfies  $k \geq \frac{n-1}{2} \geq 1$ and therefore the sum does not contain positive powers of $\lambda.$

By the construction, all the derivatives $\rho_{c}^{(s)}(c)=0, s \geq 0 .$  Also  $\rho_{c}^{(s)}(0)=0$ for all $s>0$ because $\rho_c(t)=1$ when $t$ is close to $0.$
Therefore, Leibnitz rule implies:
$$\frac{d^k}{dt^k} \big(A_M (\xi,t)\rho_c(t)\big) \big |^{t=c}_{t=0}=-\frac{d^k}{dt^k}A_M(\xi, t)\big |_{t=0}.$$
Recall that the function $A_M(\xi,t)$ is a polynomial with respect to $t.$ Therefore, if $N_0$ is its degree then all the terms in (\ref{E:sum}) with $k >N_0$ are zero and we finally obtain that $F_{\xi,1}(\lambda)$ represents in the form:

\begin{equation}\label{E:F1}
F_{\xi,1}(\lambda)=e^{i\lambda h_M(\xi)} \sum_{k=0}^{N_0} \frac{b_k(\xi)}{ \lambda^{k}}  +o\Big(\frac{1}{\lambda^{N}}\Big), \lambda \to +\infty,
\end{equation}
where $b_k(\xi)$ are the coefficients and $N$ can be taken any natural number with $N > N_0=deg_t A(\xi,t).$

Thus, we have shown that $F_{\xi,1}(\lambda)$ has a finite stationary phase expansion. Next we will prove that the remaining integrals
$F_{\xi,2}(\lambda)$ and $F_{\xi,3}(\lambda)$ vanish as $\lambda \to \infty$ to infinite order and hence contribute only to the rapidly decaying remainder.

\subsubsection{Asymptotic  of $F_{\xi,2}(\lambda)$}

It follows from the definition (\ref{E:rho}) of the function $\rho_{\xi,c}$ that
$$\Delta \rho_{\xi,c}(x)=\rho_c^{\prime\prime}(t).$$
Then, like in \ref{S:Fxi1},  we have  from Fubini theorem:
$$
F_2(\lambda)= e^{i\lambda h_M(\xi)} \int_0^c e^{i\lambda t} \rho_c^{\prime\prime}(t) A_M(\xi,t)(t)dt.
$$
The function  $\rho_c^{\prime\prime}(t)A_M(\xi,t) $ is $C^{\infty}$ and has the support inside the interval $(-c,c)$ and also $\rho_c^{\prime\prime}(t)=0$ when $t$ is close to $0$ since $\rho_c$ is constant in a neighborhood of $t=0.$ Therefore the integral $F_{\xi,2}(\lambda)$ can be viewed as the Fourier transform of a compactly supported function from $C^{\infty}(\mathbb R)$ and hence is a rapidly decaying function of the parameter $\lambda:$
\begin{equation}\label{E:F2}
F_{\xi,2}(\lambda)= o(\lambda^{-N}), \lambda \to \infty,
\end{equation}
for all $N.$

\subsubsection{ Asymptotic  of $F_3(\lambda)$}

Now turn to  the third integral in the decomposition of $I_{\xi,c}:$
$$F_{\xi,3}(\lambda)=\int\limits_{\partial D_{\xi,c}}  e^{i\lambda \langle \xi, x \rangle}  \frac{\partial}{\partial \nu_x} \rho_{\xi,c} (x) dS(x).$$

As in the case of $F_2(\lambda),$ the integration in $F_3(\lambda)$
is performed  over the relatively compact submanifold $M_{\xi,c}$ of $M.$

By the construction, $\rho_c(t)=1$ when $t$ is sufficiently close to $0.$
When $x$ is close to $a=\gamma^{-1}(\xi)$ then $t=\langle x, \xi \rangle-h_M(\xi)= \langle x-a, \xi \rangle$ is close to $0.$
Therefore
$\rho_{\xi,c}(x)=\rho_c (\langle x,\xi \rangle -h_M(\xi))=1$ when $x$ is sufficiently close to $a$

The point $a=\gamma^{-1}(\xi)$ is the only critical point of the phase function of the oscillating integral $F_3(\lambda).$
Since $\rho_{\xi,c}(x)=1$ in a neighborhood of $a$, the normal derivative of $\rho_{\xi,c}$ vanishes near the critical point $a.$
Then the stationary phase method (see \cite{Wong}) yields:
\begin{equation}\label{E:F3}
F_3(\lambda)=o(\lambda^{-N}), \lambda \to \infty,
\end{equation}
for any natural $N.$

\subsubsection{End of the proof of Lemma \ref{L:oscill}}
Now the expansion (\ref{E:decomp}) in  Lemma \ref{L:oscill} follows  from  the decomposition
$I_{\xi,c}=F_{\xi,1}+F_{\xi,2}+ F_{\xi,3}$ (\ref{E:IF123}) and asymptotic expansions (\ref{E:F1}),(\ref{E:F2}),(\ref{E:F3}). Indeed, $F_1{\xi,c}(\lambda)$  admits finite stationary phase expansion, while the second and the third ones contribute only to the rapidly decaying remainder.

In all  asymptotic relations (\ref{E:F1}),(\ref{E:F2}),(\ref{E:F3}), the remainder $R_N=o(\lambda^{-N})$ depends on the direction vector $\omega,$ i.e., $R_N=R_N(\xi,\lambda).$
However, the proof of the stationary phase expansions is based on the integration of the oscillating integrals by parts and it can be seen from there that the remainder is $o(\lambda^{-N})$ uniformly in $\xi$ with all its derivatives with respect to $\xi.$ In other words, for the remainder $R_{\xi}(\lambda)$ in (\ref{E:decomp}) holds:
$$
\lambda^N sup_{\omega} D_{\xi}^{\alpha} R(\xi,\lambda) \to 0, \lambda \to \infty.
$$
Here $D_{\xi}^{\alpha}$ is any partial derivative in $\xi$ and $N$ is an arbitrary natural number. It follows that the above asymptotic relations admit differentiation with respect to the parameter $\xi.$
Lemma is proved.

For the sake of brevity, we will use the notation
$$R_{\xi}(\lambda)=o \big( \frac{1}{\lambda^{\infty}} \big), \lambda \to \infty$$
for the remainder in (\ref{E:decomp}), if it is $o(\lambda^{-N})$ for any $N,$ uniformly with respect to $\xi$ with all the $\xi-$ derivatives.

\subsection{Differentiation the asymptotic series with respect to the direction vector}

The next step consists of the differentiation the asymptotic series (\ref{E:decomp}) with respect to the parameter $\xi.$
 To make the parameter free of the normalization, we pass in \ref{S:Renorm} to the new parameter $\omega$ with $\xi=\frac{\omega}{|\omega|}$  and replace the oscillating integral $I_{\xi,c}(\lambda)$ by a re-parametrized integral $J_{\omega,c}(\mu).$
Then in \ref{S:differ} we apply to  $J_{\omega,c}(\mu)$ the $2k$ order differential operators $\Box_{\omega}^k,$ where $\Box_{\omega}$ is the wave operator in $\omega.$
We show that after applying $\Box_{\omega}^k$ the leading term of the stationary phase expansions is obtained by multiplying the integrand
by the polynomials $T_k(x)$ which is the symbol of $\Box_{\omega}^k.$  The differentiated  integral $\Box_{\omega}^k J_{\omega,c}$ preserves the property of  having the finite stationary phase expansion and the number of terms in the expansion increases by at most $2k.$

\subsubsection{Renormalization} \label{S:Renorm}
Now awe are going to differentiate expansion (\ref{E:decomp}) with respect to the free parameter $\xi \in S^{n-1}.$
It is more convenient to deal with differentiation in the space $\mathbb R^n$ rather than on the sphere $S^{n-1},$ hence we pass from the parameter $\xi \in S^{n-1},$ to a new parameter $\omega \in \mathbb R^n:$
$$\xi=\frac{\omega}{|\omega|}, \omega \in \mathbb R^n \setminus 0.$$
In our considerations, the vector $\omega$ will be taken in an open neighborhood of the unit sphere $S^{n-1}.$

Correspondingly, we replace the large parameter
$\lambda$ in oscillating integral $I_{\xi,c}(\lambda)$ (\ref{E:I}) by
$$\mu=\frac{\lambda}{|\omega|}.$$

Then we obtain the oscillating integral
$$J_{\omega,c}(\mu):=I_{\frac{\omega}{|\omega|},c}(\mu |\omega|),$$
depending on the new parameters $\omega, \mu.$

The representation $(\ref{E:I})$ of $I_{\xi,c}$ now reads  as:
\begin{equation}\label{E:IJ}
J_{\omega,c}(\mu)=\int\limits_{M} \frac{\partial}{\partial \nu_x} e^{i \mu \langle x,\omega\rangle} \rho_{\omega,c}(x) dS(x),
\end{equation}
where
$$\rho_{\omega,c}(x)=\rho_c(\langle x-a, \frac{\omega}{|\omega|}).$$

Decomposition (\ref{E:decomp}) of $I_{\xi,c}(\lambda)$ can be rewritten, in the new parameters, as
\begin{equation} \label{E:JJ}
J_{\omega,c}(\mu)=e^{i\mu \widetilde h_M(\omega)} \widetilde Q_{\omega}\Big(\frac{1}{\mu}\Big)+\widetilde R_{\omega}(\mu)
\end{equation}
where
$$\widetilde h_M(\omega)=|\omega|h_M \Big(\frac{\omega}{|\omega|}\Big),$$
and
$\widetilde Q_{\omega}$ is a polynomial in $\frac{1}{\mu},$ with coefficients depending on $\omega:$
$$Q_{\omega}\Big(\frac{1}{\mu}\Big)=\sum\limits_{j=0}^{N_0} \frac {B_j(\omega)}{\mu^j},$$
where the leading coefficient $B_N(\omega)$ is not identically zero function.
The remainder in (\ref{E:JJ}) is a fast decaying function of $\mu^{-1}:$
$$\widetilde R_{\omega}(\mu)=o\Big(\frac{1}{ \mu^{\infty} } \Big), \mu \to \infty.$$

Notice that the index of summation in $\widetilde Q_{\omega}$ satisfies $j \geq \frac{n-1}{2} \geq 1.$

\subsubsection {Differentiation asymptotic expansion (\ref{E:JJ}) with respect to $\omega$ } \label{S:differ}
Denote
$$D_{\omega}^{\alpha}=\frac{\partial^{|\alpha|}}{\partial \omega_1^{\alpha_1} \cdots \partial \omega_n ^{\alpha_n} },  x^{\alpha}=x_1^{\alpha_1} \cdots x_n^{\alpha_n},$$
where $\alpha=(\alpha_1, ...,\alpha_n)-$ multi-index and $|\alpha|=\alpha_1 + \cdots +\alpha_n.$

\begin{lemma}\label{L:L1}
\begin{equation}\label{E:D1}
 D_{\omega}^{\alpha} J_{\omega,c}(\mu)= (i\mu)^{|\alpha|} \int\limits_{M} \frac{\partial}{\partial \nu_x}
 \Big ( e^{i \mu \langle x,\omega \rangle} x^{\alpha} \Big)  \rho_{\omega,c}(x) dS(x) + o\Big(\frac{1}{\mu^{\infty}}\Big) , \mu \to \infty.
\end{equation}
\end{lemma}
\pf

The differentiation  (\ref{E:IJ}) with respect to $\omega$ yields:

\begin{equation}\label{E:1}
D_{\omega}^{\alpha} J_{\omega,c}(\mu)=\int\limits_{M} \Big [ \frac{\partial}{\partial \nu_x} \Big( e^{i \mu \langle x,\omega \rangle}  (i\mu x)^{\alpha} \Big) \rho_{\omega,c} (x) + \mbox{other terms}\Big] dS(x),
\end{equation}
where the "other terms" contain the derivatives of $\rho_{\omega,c} (x)$ with respect to the variables $\omega_j.$

However, the function $\rho_{\omega,c}(x)$ is constant in a neighborhood of the only critical point $a=\gamma^{-1}(\frac{\omega}{|\omega|}$
of the phase function and therefore, by the stationary phase method, the terms containing the derivatives of $\rho_{\omega,c}$ contribute only to the remainder $o\big(\frac{1}{\mu^{\infty}}\big), \mu \to \infty.$


Then we obtain (\ref{E:D1}).  Lemma is proved.

Thus, we have understood derivatives of the oscillating integral. Now we turn to differentiation its  asymptotic expansion.
Direct differentiation both sides of (\ref{E:JJ}) yields
\begin{lemma}\label{L:L2}  The derivatives of $J_{\omega,c}(\mu)$ can be represented in the following form:
\begin{equation}\label{E:D2}
D_{\omega}^{\alpha}  J_{\omega,c}(\mu)=(i\mu)^{|\alpha|}  e^{i\mu \widetilde h_M(\omega)} \widetilde Q_{\omega}^{(\alpha)}
\big(\frac{1}{\mu}\big)+ o\big(\frac{1}{\mu^{\infty}}\big),
\end{equation}
where $\widetilde Q^{(\alpha)}_{\omega}$ is a polynomial with coefficients depending on $\omega$ and $\deg \widetilde Q^{(\alpha)}_{\omega} \leq N_0+|\alpha|.$
\end{lemma}

Since the left hand sides in (\ref{E:D1}) and (\ref{E:D2}) are the same, we have by equating the right hand sides:
\begin{lemma}\label{L:L3}
The following asymptotic expansion holds:
\begin{equation}\label{E:1+2}
\int\limits_M \frac{\partial}{\partial \nu_x}
\Big[ e^{i \mu \left\langle x,\omega \right\rangle} x^{\alpha} \Big]  \rho_{\omega,c}(x) dS(x)=
e^{i\mu \widetilde h_M(\omega)} \widetilde Q^{(\alpha)}_{\omega} \Big(\frac{1}{\mu} \Big)+ o\Big( \frac{1}{\mu^{\infty}} \Big),
\end{equation}
where $\widetilde Q^{(\alpha)}_{\omega}$ is a polynomial of $\deg \widetilde Q^{(\alpha)} \leq N_0+|\alpha|.$
\end{lemma}

\subsubsection{Polynomials $T_k$ and oscillating integrals $J^{k}_{\omega,c}$}
Consider the sequence of the polynomials
\begin{equation}\label{e:Tkdef}
T_k(x)=\big (x_n- \sum\limits_{j=1}^{n-1} x_j^2 \big )^k
\end{equation}
and define the oscillating integral $J^{k}_{\omega,c}(\mu)$ obtained from $J_{\omega,c}(\mu)$ by inserting the factor $T_k(x)$ under the sign of the integral:
$$
\begin{aligned}
&J^{(k)}_{\omega,c}(\mu):=\int\limits_{M} \frac{\partial}{\partial \nu_x} \Big[ e^{i \mu \left\langle x,\omega \right\rangle} T_k(x)\Big] \rho_{\omega,c}(x) dS(x) \\
&=\int\limits_{M} e^{i \mu \left\langle x,\omega \right\rangle} \Big[i\mu \left\langle \nu_x,\omega \right\rangle T_k(x) + \left\langle \nabla T_k(x), \nu_x \right\rangle \Big ]\rho_{\omega,c}(x) dS(x).
\end{aligned}
$$

\begin{corollary}  \label{C:Jexpansion}
Then the following asymptotic expansion holds:
\begin{equation}\label{E:J^k}
J^{(k)}_{\omega,c}(\mu)=e^{i\mu \widetilde h_M(\omega)} \widetilde Q_{\omega}^{(2k)} \Big(\frac{1}{\mu} \Big)+ o \Big(\frac{1}{\mu^{\infty}} \Big), \mu \to \infty,
\end{equation}
where $Q_{\omega}^{(2k)} $ is a polynomial with coefficients depending on $\omega$ and $\deg \widetilde Q_{\omega}^{(2k)} \leq N_0+2k.$
\end{corollary}
Follows from Lemma \ref{L:L3}, (\ref{E:1+2}), since $T_k(x)$ is a linear combination of $x^{\alpha}$ with $|\alpha| \leq 2k.$

\subsection {  Asymptotic of $J^{(k)}_{\omega,c}(\mu)$ as $ \mu \to \infty$}

We proceed as follows. In the previous section we have studied the transformation of the finite expansion of $J_{\omega,c}(\mu), \mu \to \infty,$  under the action of the differential operators $\Box_{\omega}^k.$ It results, on one hand, in the appearance of the factor $T_k(x)$  under the sign of integral, and, on the other hand, in enlarging the length of the expansion by at most $2k$ terms.

Next we are going to explore the dynamics of the first term of the expansion, corresponding to the minimal power of $\frac{1}{\mu}.$
Using Morse coordinates on $M$ and the osculating paraboloid condition at the distinguished point $a^0,$ we show, by establishing the order of zero of the weight function at the critical point of the phase function,  that when one applies iterates of the operator $\Box_{\omega}$ at the  point $a^0$ then the power of $\frac{1}{\mu}$ in the first  term increases faster that of the last term.

This results in disappearing the main part of the expansion after applying the operator $\Box_{\omega}$ finitely many times. In turn, it implies  vanishing on $M$  the symbol $x_n-\sum_{j=1}^{n-1}x_j^2$ of the wave operator $\Box_{\omega},$ which means that $M$ is an elliptic paraboloid  and completes the proof of Theorem \ref{T:main} .

\subsubsection{Preparations}

Let $a^0 \in M$ be the point in the formulation of Theorem \ref{T:main}, i.e., a non-degenerate point such that the osculating paraboloid at $a^0$ contacts $M$  to the order greater than four.

After applying a suitable translation and rotation, we can assume that
$$a^0=0,$$
the tangent hyperplane
$$T_{a^0}(M)=\{x_n=0 \}$$
and $$M \subset \{ x_n \geq 0\}.$$
Now the normal vector $\omega_0=\nu_{a^0}$ at the point $a^0$ becomes
$$\omega_0=(0,...,0,1).$$
We have also
$$h_M(\omega_0)=0.$$

The hypersurface $M$ can be represented in a neighborhood of $a^0=0$ as the graph
$$ x_{n}=f(x_1,\cdots, x_{n-1})$$
of a real analytic function $f.$ We have $df_0=0$ because the tangent hyperplane at $0$ is the coordinate plane $x_n=0.$ Also the third and fourth differentials ar $0$ are identically zero functions due to condition for the osculating paraboloid. Thus,
$$
f(0)=df_0=\ df^3_0=df^4_0=0.
$$

Since $a^0=0$ is a non-degenerate elliptical point of $M,$  the second differential can be reduced,
by applying a diffeomorphism of a neighborhood of 0, to the sum of squares
$$d^2 f_0(h)=h_1^2+\cdots h_{n-1}^2.$$
Then we have in a neighborhood of $0:$
\begin{equation}\label{E:f}
f(x_1,\dots,x_{n-1})=\sum_{j=1}^{n-1}x_j^2 + o(|x^{\prime}|^4).
\end{equation}

Now set
$$\omega=\omega_0=\nu_{a^0}=(0, \dots,0,1)$$
and evaluate all the above constructed objects at $\omega_0.$

Denote
$$J^{(k)}_0(\mu):=J^{(k)}_{\omega_0,c}(\mu).$$
Also we have
$$h_M(\omega_0)=0, \rho_{\omega_0, _c}(x)=\rho_c( x_n).$$

Let in expansion (\ref{E:J^k}):
$$a=a^0=0, \omega=\omega_0.$$
Then (\ref{E:J^k}) takes the following form:
\begin{equation} \label{E:J^k_0}
\begin{aligned}
&J^{(k)}_0(\mu):=\int\limits_{M} e^{i \mu x_n} \Big[ i\mu (\nu_x)_n T_k(x) + \langle \nabla T_k(x), \nu_x \rangle \Big]\rho_c(x_n) dS(x) \\
&=\widetilde Q^{(2k)}_{\omega}\Big(\frac{1}{\mu}\Big)+ o\Big(\frac{1}{\mu^{\infty}}\Big), \mu \to \infty.
\end{aligned}
\end{equation}
where $\deg \widetilde Q_{\omega}^{(2k)} \leq N_0+2k.$

\subsubsection{Oscillating integral $J^{(k)}_0(\mu)$ in local Morse coordinates on $M$ }

In this section we will express oscillating integral $J^{(k)}_0(\mu)$ in local Morse coordinates near the point $a^0.$
Then the phase function becomes quadratic and the coefficients of the stationary phase expansion express in terms of iterated Laplace operators applied to the weight function, which allows to analyse the expansion in the needed details.

First of all, we take the parameter $c$ so small, that the representation $x_n=f(x^{\prime})$ holds
on the manifold $M_{\omega_0,c} \subset M$ of integration in the integral $J^{(k)}_0(\mu).$

According to the condition of high order contact of osculating paraboloid at $a^0=0,$
we have from (\ref{E:f}) that $f$ can be represented as follows:
\begin{equation}\label{E:f}
f(x^{\prime})=|x^{\prime}|^2 + H_m(x^{\prime})+o(|x^{\prime}|^m), x^{\prime} \to 0,
\end{equation}
where $x^{\prime}=(x_1,...,x_{n-1})$ and
$$H_m=\frac{1}{m!} d^m f_0$$
is either zero or a homogeneous polynomial of degree $m >4.$

Our final goal is to prove that, in fact, $H_m=0$ for any $ m>4.$ That immediately implies that $f$ is a quadratic polynomial and, correspondingly, its graph $M$  is an elliptic paraboloid which completes the proof of Theorem \ref{T:main}.

Apply Morse lemma which says that the function $f$ can be made sum of squares by a suitable change of variables $x^{\prime}=X^{\prime}(u).$
More precisely, there exists is a (real analytic) diffeomorphism
$$X^{\prime}: B^{n-1}_{\delta} \to W_0 \subset \mathbb R^{n-1}$$
of a small ball in $\mathbb R^{n-1}$ centered at $0$  onto a neighborhood $W_0$ of $0$ in $\mathbb R^{n-1},$   such that $X^{\prime}(0)=0$ and
\begin{equation}\label{E:Morse}
(f \circ X^{\prime}) (u)=|u|^2=\sum_{j=1}^{n-1} u_j^2, u \in B^{n-1}_{\delta}.
\end{equation}
Diffeomorphism $X^{\prime}$ can be normalized
by the condition for the first differential:
$$ dX^{\prime}_0 =id,$$
so that
$$X^{\prime}(u)=u+o(u), u \to 0.$$

Also, the parameter $c$ can be taken so small that the part of the manifold $M$ on which the integration in $J^{(k)}_0(\mu)$ is performed, is contained in the graph $\{(x^{\prime},f(x^{\prime})): x^{\prime} \in W_0 \}$ over the neighborhood $W_0.$

\begin{lemma}\label{L:Phi}
$|X^{\prime}(u)|^2=|u|^2-H_m(u)+o(|u|^m), u \to 0 ,$
where $H_m$ is the homogeneous polynomial from (\ref{E:f}).
\end{lemma}
\pf
Substituting $x^{\prime}=X^{\prime}(u)$ into (\ref{E:f}) yields:
$$|u|^2=f(X^{\prime}(u))=|X^{\prime}(u)|^2 + H_m (X^{\prime}(u))+o(|X^{\prime}(u)|^m), u \to 0.$$
Since $X^{\prime}(u)=u+o(u),$ then
$$H_m(X^{\prime}(u))=H_m(u+o(u))=H_m(u)+o(|u|^m)$$ and the remainder is
$$o(|X^{\prime}(u)|^m)=o(|u+o(u)|^m)=H_m(u)+o(|u|^m).$$
Then we have
$$|u|^2=|X^{\prime}(u)|^2+ H_m(u)+o(|u|^m), u \to 0.$$
Lemma is proved.

Denote
$$X(u)=(X^{\prime}(u),|u|^2).$$
Let $J(u)$ be the Jacobian, corresponding to the change of variables $x=X(u)$ in the surface measure $dS(x):$
$$dS(x)=|J(u)| dV_{n-1}(u).$$
Now the same change of variables
in (\ref{C:Jexpansion}) leads to
\begin{equation} \label{E:Ju}
\begin{aligned}
&\int\limits_{B_{\delta}^{n-1} } e^{i \mu |u|^2} \Big ( i\mu \gamma_n\big(X(u)\big)  T_k\big(X(u)\big) + \left\langle \nabla T_k\big(X(u)\big), \gamma(X(u)) \right\rangle \Big)
\rho_c( |u|^2) |J(u)|dV_{n-1}(u) \\
&=\widetilde Q^{(2k)}_{\omega}\Big(\frac{1}{\mu}\Big)+ o\Big(\frac{1}{\mu^{\infty}}\Big), \mu \to \infty.
\end{aligned}
\end{equation}

\subsubsection{ Asymptotic of the weight function in $J_0^k$ near $u=0$}

In order to understand the asymptotic of $J^k_0(\mu)$ as $\mu \to \infty,$
let us explore the behavior of functions under the sign of the integral $J^k_0(\mu)$ near $u=0.$

We have due to (\ref{L:Phi}):
\begin{equation}\label{E:TkX}
T_k(X(u))=(|u|^2- |X^{\prime}(u)|^2)^k=H^k_m(u)+o(|u|^{km}|).
\end{equation}
Remind also that in the new coordinates the distinguished point is  $a^0=0.$
Also $X(0)=0$ and
$$\gamma(X(0))=\gamma(a^0)=(0,...,0,1),$$
\begin{equation} \label{E:gamma}
\begin{aligned}
&\gamma_j(u)=l_j(u)+o(|u|), j=1,...,n-1,  \\
&\gamma_n(X(u))=1+o(1), u \to 0,
\end{aligned}
\end{equation}
where $l_j(u)$ are  linear functions. Remind that $\gamma(x)=(\gamma_1(x), ...,\gamma_n(x))=\nu_x$ is Gaussian mapping on $M.$

Then we have from \ref{E:TkX}:
\begin{equation}\label{E:T_k}
\begin{aligned}
&\gamma_n(X(u))T_k(X(u))=(1+o(1))(|u|^2- |X^{\prime}(u)|^2)^k
&\\=H_m^k(u)+o(|u|^{km}) \\
\end{aligned}
\end{equation}
From definition (\ref{E:Tkdef}) of the polynomials $T_k:$
$$\nabla T_k(x)=k(x_n- |x^{\prime}|^2)^{k-1} (-2x^{\prime}, 1).$$
Since $X^{\prime}(u)=u+o(u)$ and $|u|^2-|X^{\prime}|^2=H_m(u)+o(|u|^{m})$  we have
$$
\begin{aligned}
&\nabla T_k(X(u))=k(|u|^2-|X^{\prime}(u)|^2)^{k-1} (-2X^{\prime}(u),1) \\
&=k (H_m^{k-1} (u)+o(|u|^{(k-1) m})) (-2u+o(|u|), 1+o(1)), u \to 0.
\end{aligned}
$$

Taking into account (\ref{E:gamma}), we compute the inner product of the vectors $\nabla T_k(X(u))$  and $\gamma(X(u),$ and obtain

\begin{equation}\label{E:nablaT_k}
\begin{aligned}
&\langle \nabla T_k(X(u)), \gamma\big(X(u)\big) \rangle \\
&= \Big (k H_m^{k-1}(u)+ o(|u|^{(k-1) m}\Big) \Big( -2\sum\limits_{j=0}^{n-1} l_j(u)u_j +o(|u|^2) +  (1+o(1))\Big) \\
&=kH_m^{k-1}(u)+ o(|u|^{(k-1) m}), u \to 0.
\end{aligned}
\end{equation}

Represent integral $J^k_0(\mu)$ in (\ref{E:J^k_0}) as the sum of the two integrals
$$J^k_0(\mu)=i \mu J_1(\mu)+J_2(\mu),$$
where
\begin{equation}\label{E:J12}
\begin{aligned}
&J_1(\mu):=\int\limits_{B_{\delta}^{n-1} } e^{i \mu |u|^2}  \gamma_n\big(X(u)\big)  T_k\big(X(u)\big) \rho_c( |u|^2) |J(u)|dV_{n-1}(u), \\
&J_2(\mu)=\int\limits_{B_{\delta}^{n-1} } e^{i \mu |u|^2} \left\langle \nabla T_k\big(X(u)\big), \gamma\big(X(u)\big) \right\rangle \rho_c( |u|^2) |J(u)|dV_{n-1}(u).
\end{aligned}
\end{equation}

Thus, substituting the expansions (\ref{E:T_k}) and (\ref{E:nablaT_k})  into (\ref{E:J12}) leads to

\begin{equation}\label{J:J1}
J_1(\mu)=
\int\limits_{B_{\delta}^{n-1} } e^{i \mu |u|^2} \big(H_m^{k}(u)+o(|u|^{k m}\big) \rho_c( |u|^2) |J(u)|dV_{n-1}(u)
\end{equation}
and
\begin{equation}\label{E:J2}
J_2(\mu)=
\int\limits_{B_{\delta}^{n-1} } e^{i \mu |u|^2} \big( kH_m^{k-1}(u)+ o(|u|^{(k-1) m})\big) \rho_c( |u|^2) |J(u)|dV_{n-1}(u).
\end{equation}

Now we are going to determine the leading terms of the asymptotic expansion,
as $\mu \to \infty$, of each oscillating integral $J_1(\mu)$ and $J_2(\mu)$ separately.

\subsubsection {The leading term of the asymptotic of $J_1(\mu), \mu \to \infty$}
Integral $J_1(\mu)$ is an oscillating integral
$$J_1(\mu)=i\mu \int\limits_{B_{\delta}^{n-1}} e^{i \mu |u|^2} p(u) dV_{n-1}(u)$$
in a neighborhood of $0 \in \mathbb R^{n-1},$ with the large parameter $\mu,$ the quadratic phase function $|u|^2$
and the weight function
$$p(u):=(H_m^{k}(u)+o(|u|^{k m}) \rho_c( |u|^2) |J(u)|.$$

According to the stationary phase method (see, e.g., \cite{Wong}), $J_1(\mu)$ expands into the asymptotic series:
$$
J_1(m)=\sum\limits_{j=0}^{\infty} c_j \frac{\big(\Delta^j p\big)(0)}{\mu^j},$$
where
$\Delta $ is the Laplace operator with respect to the variable $u$ and $c_j$ are certain nonzero numerical coefficients.

Since $H_m^k$ is a homogeneous polynomial of degree $mk,$ the function $\Delta^j H_m^k$ is a homogeneous polynomial of degree
$$\deg \Delta^j H_m^k=mk-2j.$$
Therefore
$$\Delta^j p(0)=0$$
for $2j<mk.$
and
\begin{equation}\label{E:J1c}
J_1(\mu)=\sum\limits\limits_{j \geq \frac{mk}{2}} c_j \frac{\big(\Delta^j p \big)(0)}{\mu^j}.
\end{equation}
with

\subsubsection {The leading term of the asymptotic of $J_2(\mu), \mu \to \infty$}

From (\ref{E:J2}):
$$
J_2(\mu)=
\int\limits_{B_{\delta}^{n-1} } e^{i \mu |u|^2} q(u) dV_{n-1}(u).
$$
where
\begin{equation} \label{E:qqq}
q(u):=\big(kH_m^{k-1}(u)+ o(|u|^{(k-1) m})\big) \rho_c( |u|^2) |J(u)|, u \to 0.
\end{equation}
The leading term of the asymptotic of $q(u)$ as $u \to 0$ is determined by the term $H_m^{k-1}(u)$ which is a homogeneous polynomial of degree $m(k-1).$
Therefore
$$
\big(\Delta^j q \big)(0)=0
$$
for $2j < m(k-1)$ and hence
\begin{equation} \label{E:J2c}
J_2(\mu)=\sum\limits\limits_{j \geq \frac{m(k-1)}{2} } c_j \frac{\big(\Delta^j q \big)(0)}{\mu^j}.
\end{equation}

\subsubsection {The leading term of the asymptotic of $J^k_0(\mu), \mu \to \infty$}

Now we are able to determine the leading term of the asymptotic expansion for the sum $J^k_0(\mu)=i\mu J_1(\mu)+J_2(\mu):$
\begin{lemma} \label{E:leadterm}  Take $k$ odd, $k=2\alpha+1.$ Then
\begin{equation}\label{E:leadterm1}
J_0^k(\mu)= C \frac{\big(\Delta^{m\alpha}H_m^{2\alpha}\big)(0) }{\mu^{m\alpha}}+ o\Big(\frac{1}{\mu^{m\alpha}}\Big), \mu \to \infty,
\end{equation}
where $C$ is a nonzero constant.
\end{lemma}
\pf
Since $J^k_0(\mu) = i\mu J_1(\mu)+J_2(\mu),$ we have from (\ref{E:J1c}), (\ref{E:J2c}):
$$J^k_0(\mu)= i\mu \sum\limits\limits_{j \geq \frac{ mk}{2}} c_j \frac{\big(\Delta^j p \big)(0)}{\mu^j}+
\sum\limits\limits_{j \geq \frac{m(k-1)}{2}} c_j \frac{\big(\Delta^j p \big)(0)}{\mu^j}..$$

Let us compare the minimal degrees of $\frac{1}{\mu}$ contributing by each sum.
The lower bound for the degree $\frac{1}{\mu}$ in  $i\mu J_1(\mu)$  is
$$j_1=\frac{mk}{2}-1 =m\alpha +\frac{m}{2}-1,$$
while the minimal degree coming from the second sum is
$$j_2=\frac{m(k-1)}{2}=m\alpha.$$
Since $m >4$ we have $j_1> m \alpha +1 > j_2$ and hence the leading term of the expansion for $J^k_0(\mu)$
comes from $J_2$ only, i.e.,
$$J_0^k(\mu)=c_{j_2} \frac{\big(\Delta^{j_2} q \big)(0)}{\mu^{j_2}}+ o\Big(\frac{1}{\mu^{j_2}} \Big), \mu \to \infty.$$
That is exactly what Lemma claims because
$$
\rho_c( |u|^2) |J(u)|=1+o(1), u \to 0.
$$
and hence from (\ref{E:qqq}), with $k=2\alpha+1,$ follows:
$$
\begin{aligned}
&\big(\Delta^{j_2} q \big)(0)=\Delta^{m\alpha} \Big( \big( (2\alpha+1)H_m^{2\alpha}(u)+ o(|u|^{2 m\alpha }) \big) \rho_c( |u|^2) |J(u)| \Big)(0) \\ &= C \big (\Delta^{m\alpha}H_m^{2\alpha}\big)(0).
\end{aligned}
$$
Lemma is proved.
\subsection{End of the proof of Theorem \ref{T:main}}
\begin{lemma}\label{L:Delta}
$$
(\Delta^{m\alpha}H_m^{2\alpha})(0) =0,
$$
for sufficiently large $\alpha.$
\end{lemma}
\pf
Compare asymptotic expansions (\ref{E:J^k_0}) and (\ref{E:leadterm1}), with respect to the large parameter $\mu$ of the same oscillating integral
\begin{equation} \label{E:expanexpan}
J^k_0(\mu)=\widetilde Q^{(2k)}\big(\frac{1}{\mu}\big)+ o\big(\frac{1}{\mu^{\infty}}\big)=C \frac{\big(\Delta^{m\alpha}H_m^{2\alpha}\big)(0) }{\mu^{m\alpha}}+ o\big(\frac{1}{\mu^{m\alpha}}\big),  \mu \to \infty.
\end{equation}
By ( \ref{E:J^k_0}) the degree of polynomial $Q^{(2k)}$ is at most $N_0+2k.$
Therefore the term in the first decomposition (\ref{E:expanexpan}), with the highest power of $\mu^{-1},$ has the form
$$ \frac{B}{\mu^N},$$
where $B$ is a constant and
the exponent in the denominator satisfies
$$N \leq N_0+2k=N_0+4\alpha.$$

On the other hand, we see from (\ref{E:expanexpan}) that the minimal degree of  $\frac{1}{\mu}$ in the expansion
is $m\alpha.$
Since $m > 4,$ we have
$$m\alpha > N_0+4\alpha,$$
as soon as
$$\alpha >\frac{N_0}{m-4}.$$
Thus, for large $\alpha,$ the  lowest degree term  $\mu^{-m\alpha}$ in the second expansion in (\ref{E:expanexpan}) has the degree higher than the highest degree term  in the first expansion of the same function $J^k_0(\mu).$

Therefore the term with $\mu^{-m\alpha}$ is not presented in (\ref{E:expanexpan}) at all , which means that the coefficient in front of that term is zero:
$$(\Delta^{m\alpha}H_m^{2\alpha})(0)=0.$$
Lemma is proved.

\begin{lemma} \label{L:Hm=0}
$H_m=0.$
\end{lemma}
\pf
Represent the variable $u \in \mathbb R^{n-1}$ in the spherical coordinates $u=r\theta, |\theta|=1.$
Then $H_m(r\theta)=r^mH_m(\theta)$ and the power of $H_m$ is
$$H_m^{2\alpha}(u)=r^{2\alpha m} H_m^{2\alpha}(\theta).$$
Consider the spherical average
$$\widetilde{H_m^{2\alpha}}(r)=\int\limits_{|\theta|=1}H_m^{2\alpha}(r\theta)dA(\theta)= A r^{2m\alpha},$$
where
\begin{equation} \label{E:A}
A:= \int\limits_{|\theta|=1} H_m^{2\alpha}(\theta)dA(\theta).
\end{equation}
Since $\Delta$ is invariant with respect to rotations, we have from Lemma \ref{L:Delta}:
$$ 0=(\Delta^{m\alpha}H_m^{2\alpha})(0) =(\Delta^{m\alpha}\widetilde{H_m^{2\alpha}})(0) = C \cdot A,$$
where $$C=\prod_{j=1}^{m\alpha}\big((2m\alpha-2j+2)(2m\alpha-2j+1)+n-2 \big) \neq 0.$$

Thus, $A=0.$ The function under the integral in (\ref{E:A}) is non-negative, therefore $H_m(\theta)=0$ identically on the unit sphere and then the homogeneous polynomial $H_m(u)=0$ identically. Lemma is proved.

Now we can complete the proof of Theorem \ref{T:main}. We have denoted by $H_m=\frac{1}{m!}d^mf_0$ the  first nonzero homogeneous polynomial of degree $m>4$  in the Taylor formula  if such a term exists and Lemma \ref{L:Hm=0} says that it does not:
$$H_m=d^m f_0=0.$$
Therefore, all the differentials
$$d^m f_0=0, m>4.$$
Since, by the assumption, we have $d^3 f_0=d^4 f_0=0$ and $f$ is real analytic, then
$$f(x^{\prime})=|x^{\prime}|^2.$$
That means that the hypersurface $M,$ which is the graph of $f,$ is an elliptic paraboloid. Theorem \ref{T:main} is proved.

\section{Concluding remarks}

The proof of Theorem \ref{T:main} is based on the equivalence of  polynomial integrability of a hypersurface and finiteness of the stationary phase expansion, with respect to a large parameter, of oscillating integrals with linear phase function (Fourier integrals).

Study of finite stationary phase expansions is of independent interest (see \cite{Bernard} and references there).
In our case, we deal with not a single phase function but with  a parametric family of the phase functions $\varphi_{\omega}(x)=\langle \omega, x \rangle,$  depending on the direction vector $\omega.$
What we essentially have proven is that the finiteness of the stationary phase expansion  for such a family of oscillating integrals on $M$ does not hold except for quadrics in odd-dimensional spaces.

The assumption of existence a distinguished point $a^0$ with osculating paraboloid of high order contact, restricts the result since it rules out ellipsoids and two-sheet hyperboloids. Getting rid of this restriction would allow to obtain complete characterization of convex quadrics in odd-dimensional spaces in terms of finite stationary phase expansions.

Bar-Ilan University and Holon Institute of Technology; Israel

{\it E-mail}: agranovs@math.biu.ac.il

\end{document}